\magnification=1200
\baselineskip=14pt plus 1pt minus 1pt
\tolerance=1000\hfuzz=1pt 


\font\bigfont=cmr10 scaled\magstep3
\font\secfont=cmr10 scaled\magstep1

\def\section#1#2{\vskip32pt plus4pt \goodbreak \noindent{\bf\secfont#1. #2}
        \xdef\currentsec{#1} \global\eqnum=0 \global\thmnum=0}
\def\appendix#1#2{\vskip32pt plus4pt \goodbreak \noindent
                {\bf\secfont{Appendix #1. }#2}
        \xdef\currentsec{#1} \global\eqnum=0 \global\thmnum=0}

\def\leftheader{}
\def\rightheader{}
\def\header{\headline={\tenrm\ifnum\pageno>1 
                        \ifodd\pageno\rightheader
                         \else\noindent\leftheader 
                          \fi\else\hfil\fi}}

\newcount\thmnum
\global\thmnum=0
\def\prop#1#2{\global\advance\thmnum by 1
        \xdef#1{Proposition \currentsec.\the\thmnum}
        \bigbreak\noindent{\bf Proposition \currentsec.\the\thmnum.}
        {\it#2} }
\def\defin#1#2{\global\advance\thmnum by 1
        \xdef#1{Definition \currentsec.\the\thmnum}
        \bigbreak\noindent{\bf Definition \currentsec.\the\thmnum.} 
        {#2} }
\def\lemma#1#2{\global\advance\thmnum by 1
        \xdef#1{Lemma \currentsec.\the\thmnum}
        \bigbreak\noindent{\bf Lemma \currentsec.\the\thmnum.} 
        {\it#2} }
\def\thm#1#2{\global\advance\thmnum by 1
        \xdef#1{Theorem \currentsec.\the\thmnum}
        \bigbreak\noindent{\bf Theorem \currentsec.\the\thmnum.} 
        {\it#2} }
\def\cor#1#2{\global\advance\thmnum by 1
        \xdef#1{Corollary \currentsec.\the\thmnum}
        \bigbreak\noindent{\bf Corollary \currentsec.\the\thmnum.} 
        {\it#2} }
\def\conj#1#2{\global\advance\thmnum by 1
        \xdef#1{Conjecture \currentsec.\the\thmnum}
        \bigbreak\noindent{\bf Conjecture \currentsec.\the\thmnum.} 
        {\it#2} }

\def\proof{\medskip\noindent{\it Proof. }}

\newcount\eqnum
\global\eqnum=0
\def\num{\global\advance\eqnum by 1
        \eqno({\rm\currentsec}.\the\eqnum)}
\def\eqalignnum{\global\advance\eqnum by 1
        ({\rm\currentsec}.\the\eqnum)}
\def\ref#1{\num  \xdef#1{(\currentsec.\the\eqnum)}}
\def\eqalignref#1{\eqalignnum  \xdef#1{(\currentsec.\the\eqnum)}}

\def\title#1{\centerline{\bf\bigfont#1}}

\newcount\subnum
\def\Alph#1{\ifcase#1\or A\or B\or C\or D\or E\or F\or G\or H\fi}
\def\subsec{\global\advance\subnum by 1 
        \vskip12pt plus4pt \goodbreak \noindent
        {\bf \currentsec.\Alph\subnum.}  }
\def\newsubsec{\global\subnum=1 \vskip6pt\noindent  
        {\bf \currentsec.\Alph\subnum.}  }
\def\today{\ifcase\month\or January\or February\or March\or 
        April\or May\or June\or July\or August\or September\or 
        October\or November\or December\fi\space\number\day, 
        \number\year}

\input amssym.def
\input amssym.tex

\def\ol{\overline}

\def\tr{\mathop{\rm tr}\nolimits}
\def\str{\mathop{\rm Str}\nolimits}

\def\re{\mathop{\rm Re}\nolimits}
\def\im{\mathop{\rm Im}\nolimits}

\def\bC{{\Bbb C}}
\def\bN{{\Bbb N}}

\def\bR{{\Bbb R}}

\def\bZ{{\Bbb Z}}

\def\cA{{\cal A}}

\def\cD{{\cal D}}
\def\cF{{\cal F}}
\def\cH{{\cal H}}

\def\cL{{\cal L}}

\def\cP{{\cal P}}

\def\fA{{\frak A}}

\def\fK{{\frak K}}

\def\fP{{\frak P}}

\def\fT{{\frak T}}

\chardef\o="1C

\header
\def\leftheader{\centerline{S. KLIMEK and A. LE\'SNIEWSKI}}
\def\rightheader{\centerline{QUANTUM FIELD THEORY}}

\def\intsec{I}
\def\kronsec{II}
\def\bosesec{III}
\def\fermsec{IV}
\def\skmssec{V}
\def\tauberapp{A}
\def\examplapp{B}

{\baselineskip=12pt
\nopagenumbers
\line{\hfill \bf HUTMP 95/445}
\line{\hfill \today}
\vfill
\title{Quantized Kronecker flows}
\vskip 1cm
\title{and almost periodic quantum field theory}
\vskip 1in
\centerline{{\bf S\l awomir Klimek}$^*$
\footnote{$^1$}{Supported in part by the National Science Foundation under
grant DMS--9500463}
and {\bf Andrzej Le\'sniewski}$^{**}$\footnote{$^2$}
{Supported in part by the National Science Foundation under grant
DMS--9424344}}
\vskip 12pt
\centerline{$^*$Department of Mathematics}
\centerline{IUPUI}
\centerline{Indianapolis, IN 46205, USA}
\vskip 12pt
\centerline{$^{**}$Lyman Laboratory of Physics}
\centerline{Harvard University}
\centerline{Cambridge, MA 02138, USA}
\vskip 1in\noindent
{\bf Abstract.} We define and study the properties of the infinite
dimensional quantized Kronecker flow. This $\bC^*$-dynamical system
arises as a quantization of the corresponding flow on an infinite
dimensional torus. We prove an ergodic theorem for a class of quantized
Kronecker flows. We also study the closely related almost periodic
quantum field theory of bosonic, fermionic and supersymmetric particles.
We prove the existence and uniqueness of KMS and super-KMS states for
the $\bC^*$-algebras of observables arising in these theories. 
\vfill\eject}

\section\intsec{Introduction}
\medskip
In this paper, we introduce and study the properties of quantized infinite
dimensional Kronecker flows. Very much like its classical counterpart, a
quantized Kronecker flow is defined in terms of an infinite sequence
$\Omega$ of frequencies satisfying certain genericity assumptions. We
define a natural $\bC^*$-algebra of observables and show that it can be
identified with an infinite tensor product of standard Toeplitz algebras.
A quantized Kronecker flow is a one parameter group of automorphisms of this
$\bC^*$-algebra.

Our work has been inspired by the recent preprint [BC]. The structures
studied in [BC] can be interpreted as an example of a quantized
Kronecker flow with the set of frequencies $\Omega$ equal
to $\{\log p:\; p\hbox{ prime number}\}$.

The infinite tensor product of standard Toeplitz $\bC^*$-algebras referred
to above arises naturally when one quantizes the infinite dimensional
Kronecker flow on a Bohr compactification of $\bR$. Even though there is
no parameter in this theory which would play the role of Planck's
constant, one can introduce a natural concept of the classical limit
[C1], [S1], [Z]. We prove that, under additional assumptions on
$\Omega$, the classical limit of the quantized Kronecker flow exists.
The proof of this theorem relies on an Ingham type tauberian theorem.
Furthermore, we study the ergodic properties of that quantum Kronecker
flow. We show that whenever the classical limit exists, the quantized
Kronecker flow is quantum ergodic in the sense of Zelditch [Z].
  
Infinite dimensional Kronecker flows lead to models of free
quantum field theory in one space dimension. In these
field theories, the field operators are almost periodic functions
of the space coordinate $x$. There is a natural notion of a mean
in the theory of almost periodic functions, the Bohr mean, which plays
the role of the integral. Using it, we carry over much of the formalism
of quantum field theory to the almost periodic setup. The construction
of an almost periodic field theory requires ``doubling" the Hilbert
space of the Kronecker flow. On this Hilbert space, we define the
field and momentum operators, which are fundamental objects in
canonical quantum field theory. It turns out that the dynamics
of the quantum Kronecker flow and the almost periodic
wave equation are essentially the same. We focus our discussion on
the theory of free fields.. Interacting (i.e. nonlinear) field theories
exhibit some new striking phenomena and will be discussed in a
future publication.

Furthermore, we extend the construction of almost periodic field theory
to incorporate fermionic and supersymmetric (i.e. $\bZ_2$-graded) fields.
The supersymmetric extension leads to a natural construction of a
Fredholm module over the algebra of observables associated with the
almost periodic quantum field and the corresponding super-KMS functional.
Super-KMS functionals [K], [JLW], appear naturally in the $\bZ_2$-graded
entire cyclic cohomology theory, but have been studied less intensively
than their nongraded antecedents, namely KMS states. We prove the 
uniqueness of the super-KMS functionals for the supersymmetric model
of an almost periodic quantum field theory.

The paper is organized as follows. In Section {\kronsec} we introduce
the notation and recall basic facts from the theory of almost
periodic functions. We then study the ergodic properties
of the quantized almost periodic Kronecker flow. The main technical
input is a variant of Ingham's tauberian theorem proved in Appendix
{\tauberapp}. As it turns out, ergodicity of the quantized Kronecker
flow depends on the growth properties of the set $\Omega$. Examples of
ergodic Kronecker systems are given in Appendix {\examplapp}.
We start Section {\bosesec} with a discussion of the canonical formalism
of the almost periodic classical field theory and illustrate it
with an analysis of the almost periodic wave equation. Then we show
how to formulate the quantum version of the almost periodic
wave equation. Free fermionic and supersymmetric models
are introduced in Section {\fermsec}. Finally,
in Section {\skmssec}, we prove the uniqueness of the super-KMS
functionals for the free supersymmetric almost periodic
quantum field theory.

\section\kronsec{Kronecker flows}
\medskip
In this section we introduce the central concept of this paper, namely
the infinite dimensional Kronecker flow. After a brief summary of the 
classical theory, we present its quantum version and study the properties
of the resulting dynamics.
\newsubsec
First, we review some facts from the theory of almost periodic
functions on $\bR$ and fix our notation.
\defin\kroneckerdefine{A countable ordered subset $\Omega$ of $\bR$ is
called a {\it Kronecker system} if it satisfies conditions 1-4 below.
\item{1.} {$0\notin\Omega$.} 
\item{2.}{Let $\Omega=\Omega_{+}\cup\Omega_-$, where $\Omega_+$ and
$\Omega_-$ are the subsets of positive and negative elements of
$\Omega$, respectively. Then $\Omega$ is even i.e. $\Omega_-=-\Omega_+$.} 
\item{3.} {Let $\omega_n$, $n=1,\ldots$, be the elements of $\Omega_+$
listed in increasing order. Then $\omega_n\to\infty$, as $n\to\infty$.} 
\item{4.} {The elements of $\Omega_+\ $ are algebraically independent
over $\bZ$.}}
\medskip
Two natural examples of a Kronecker system arise as follows.

\noindent
{\it Example 1.} Let $\fK$ be an algebraic number field and $\fP$
the the set of its prime ideals. Set
$$
\Omega=\{\pm\log NP:\ P\in\fP\},
$$
where $NP$ denotes the norm of $P$. Then $\Omega$ is a Kronecker
system. In particular, the unique factorization property of ideals implies 
that $\Omega_+\ $ consists of numbers algebraically independent over 
$\bZ$. This example is taken from [BC], [C2], [J], [S2].

\noindent
{\it Example 2.} Let $m$ be a transcendental real number, and let
$$
\Omega=\{\pm\sqrt{n^2+m^2}:\;n\in\bN\}.
$$
Then $\Omega$ is a Kronecker system. This example is motivated
by two dimensional quantum field theory.

Let $\ol\Omega\subset\bR$ be the set of linear combinations of elements
of $\Omega$ with coefficients in $\bZ$. In other words, $\ol\Omega$ is
the free abelian group generated by $\Omega_+$:
$$
\ol\Omega=\bZ\left[\Omega_+\right]\ .
$$ 
We will use the notation $\cA\cP(\Omega)$ for the vector space of
continuous almost periodic functions on $\bR$ with frequencies in
$\Omega$. This means that a function $f\in \cA\cP(\Omega)$ has the
following uniformly convergent Fourier expansion:
$$
f(x)=\sum_{\omega\in\Omega} f_\omega e^{i\omega x}.
$$
Likewise, $\cA\cP(\ol\Omega)$ will denote the space of continuous
almost periodic functions on $\bR$ with frequencies in $\ol\Omega$.
Unlike $\cA\cP(\Omega)$, the space $\cA\cP(\ol\Omega)$ forms a unital
commutative $\bC^*$-algebra. This $\bC^*$-algebra can be identified
with the $\bC^*$-algebra of continuous functions on the following Bohr 
compactification of $\bR$. 

Let $\ol\bR_\Omega$ denote the infinite cartesian product of unit circles,
$\ol\bR_\Omega=\prod_{\omega\in\Omega_+}S^1$, equipped with the Tikhonov
topology. The embedding
$$
\bR\ni x\longrightarrow\prod_{\omega\in\Omega_+}\exp{i\omega x}
\in\prod_{\omega\in\Omega_+}S^1
\ref{\embedref}
$$
induces an isomorphism $C(\ol\bR_\Omega)\simeq \cA\cP(\ol\Omega)$,
see [HR].

The product of Lebesgue measures on $S^1$ defines an integral
$\int_{ap}$ on $C(\ol\bR_\Omega)$ and, consequently, on
$\cA\cP(\ol\Omega)$. Explicitly, in terms of Fourier series we have
$$
\int_{ap}\sum_{\eta\in\ol\Omega}f_\eta e^{i\eta x}\,dx=
f_0\ .\ref{\apintref}
$$

For later reference, we note that the almost periodic Dirac's delta
distribution $\delta_\Omega$ defined by
$$
\delta_\Omega(x)=\sum_{\omega\in\Omega}e^{i\omega x}
\ref{\apdeltaref}
$$
is the Schwartz kernel of the projection
$$
\cA\cP(\ol\Omega)\ni f\longrightarrow\int_{ap}\delta_\Omega(x-y)f(y)\,dy
\in \cA\cP(\Omega)\ .
$$
This projections ``forgets'' all terms in the Fourier series whose
frequencies are not in $\Omega$.

The embedding \embedref\ defines a Kronecker type flow $\alpha_t$ on
$\ol\bR_\Omega$ given by
$$
\alpha_t\big(\prod_{\omega\in\Omega_+}e^{ix_\omega}\big)=
\prod_{\omega\in\Omega_+}e^{ix_\omega+it\omega}\ .\ref{\flowref}
$$ 
As a consequence of our assumptions on $\Omega$, this flow is ergodic.

\subsec
We will now construct a quantization of the classical dynamical system
$(\cA\cP(\ol\Omega),\alpha_t)$, and study the ergodic properties of
the resulting quantum Kronecker flow. The quantization will be given
in terms of an algebra of operators on a Hilbert space, the ``algebra
of observables''.

Set $\cH_+=l^2(\Omega_+)$, and consider the bosonic Fock space
$\cF_b\cH_+$ defined as the symmetric tensor algebra over $\cH_+$,
$$
\cF_b\cH_+=\bigoplus_{n=0}^{\infty}S^n\cH_+\ ,
$$
where $S^n\cH_+$ is the $n^{\hbox{th}}$ symmetric tensor power of $\cH_+$
(with $S^0\cH_+=\bC$). The vector $v_0=(1,0,0,\ldots)\in \cF_b\cH_+$ is
called the {\it vacuum}. The Hilbert space $\cF_b\cH_+$ can be naturally
identified with $l^2(\bN[\Omega_+])$, where $\bN[\Omega_+]$ is
the nonnegative cone in the lattice $\bZ[\Omega_+]$. 
In the example of an algebraic number field $\fK$, the set
$\bN[\Omega_+]$ can be identified with the set of all ideals of $\fK$.

Alternatively, there is a natural isomorphism of 
$\cF_b\cH_+$ with an infinite tensor product
$$
\cF_b\cH_+\simeq\bigotimes_{\omega\in\Omega_+}l^2(\bN[\omega])\ .
\ref{\kronprodref}
$$
In von Neumann's terminology, if $e_n(\omega),\ n=0,1,2,\ldots$, is the
canonical basis in $l^2(\bN[\omega])$, then the above tensor product is
the incomplete tensor product associated with the sequence of vectors
$(e_0(\omega),e_0(\omega),e_0(\omega),\ldots)$. Furthermore, the Fourier
transform allows us to identify the space $l^2(\bZ[\Omega_+])$ with
$L^2(\ol\bR_\Omega)$, and the Fock space  $\cF_b\cH_+\simeq l^2
(\bN[\Omega_+])$ with the closed subspace $L^2_+(\ol\bR_\Omega)$ of
$L^2(\ol\bR_\Omega)$ consisting of functions with nonnegative frequencies.

Let $P$ be the orthogonal projection onto $L^2_+(\ol\bR_\Omega)
\subset L^2(\ol\bR_\Omega)$. Every $f\in C(\ol\bR_\Omega)$ defines a
{\it Toeplitz operator} $T(f)$ on $L^2_+(\ol\bR_\Omega)\simeq\cF_b\cH_+$
by
$$
T(f)=PM(f)P,
$$
where $M(f)$ is the operator on $L^2(\ol\bR_\Omega)$ of multiplication
by $f$. Recall that $T(f)$ is continuous in $f$, $\Vert T(f)\Vert
\leq\Vert f\Vert_\infty$, where $\Vert f\Vert_\infty$ denotes the
sup norm of $f$.  

Let $\cA_+$ be the $C^*$-algebra generated by the Toeplitz operators.
It is not difficult to see that $\cA_+$ coincides with the (reduced)
$\bC^*$-algebra of the semigroup $\bN[\Omega_+]$.
For $\eta\in\bN[\Omega_+]$, let $e(\eta)$ denote the canonical basis
element in $l^2(\bN[\Omega_+])\simeq \cF_b\cH_+$. Let $H_+$ be an
unbounded, self-adjoint operator in $\cF_b\cH_+$ defined by
$$
H_+e(\eta)=\eta e(\eta),
\ref{\hplusdefref}
$$
and let $U(t)=e^{itH_+}$ be the corresponding one parameter group of
unitary operators. The pair $(\cA_+,U(t))$ is a quantization
of $(C(\ol\bR_\Omega),\alpha_t)$ which we call the {\it quantum
Kronecker flow}. 

Recall that the standard Toeplitz $\bC^*$-algebra $\fT$ is the
$\bC^*$-algebra generated by a single generator $u$ satisfying the
relation $u^*u=I$. The following proposition can be proved by the
method used in the proofs of Propositions 7 and 8 in [BC].

\prop\toeplitzprop{
\item{1.} The $\bC^*$-algebra $\cA_+$ is an infinite tensor
product of Toeplitz $\bC^*$-algebras $\fT_\omega$,
$$
\cA_+=\bigotimes_{\omega\in\Omega_+}\fT_\omega\ ,
$$
where $\fT_\omega$ is generated by the unilateral shift
$u_\omega=T(e^{ix_\omega})$.
\item{2.} For every $\beta>0$ there exists a unique KMS$_\beta$
state for $(\cA_+,U(t))$.
}
\vskip 5mm
\subsec
We will show now that the quantum dynamical system constructed
above is in fact a quantization of the classical Kronecker flow.
Even though there is no Planck's constant in this theory, one can
still introduce a natural concept of its {\it classical limit}.
Such a construction of the classical limit was originally proposed
in [C1], [S1], [Z], and consists in the following. If $a$ is an
operator on $\cF_b\cH_+$ and $E>0$, we set
$$
\tau_E(x):={1\over N(E)}\sum_{\bN[\Omega_+]\ni\eta\leq E}
(e(\eta),ae(\eta)),
\ref{\tauedefref}
$$
where $N(E)$ is the number of eigenvalues of $H_+$ less than or equal to
$E$. The theorem which we are about to formulate describes the classical
limit of $(\cA_+,U(t))$.

We will need additional assumptions on the set $\Omega_+$ to guarantee
the existence of the classical limit. Specifically,
assume that for every $s>0$,
$$
\theta(s):=\sum_{n=1}^\infty e^{-s\omega_n}<\infty .\ref{\thetasumref}
$$
This implies that the following $\zeta$-type function
$$
\zeta_\Omega(s):=\prod_{n=1}^\infty(1-e^{-s\omega_n})^{-1}
\ref{\zetadefref} 
$$
converges for all $s>0$. Expanding each term of $\zeta_\Omega(s)$
in a power series and multiplying out the terms, we can express
$\zeta_\Omega(s)$ as the following Lebesque-Stietljes integral:
$$
\zeta_\Omega(s)=1+\int_0^\infty e^{-sx}\,dN(x) ,
$$
where, as above, $N(x)$ is the counting function for the
eigenvalues of $H_+$. Equivalently, we can write this formula as
$$
e^{\phi(s)}=1+\int_0^\infty e^{-sx}\,dN(x) ,
\ref{\laplacetransref}
$$
where $\phi(s):=-\sum_{n=1}^\infty\log(1-e^{-s\omega_n})$.
In Appendix {\tauberapp} we study {\laplacetransref} recast in the
following form:
$$
{e^{\phi(s)}\over s}=\int_0^\infty e^{-sx}(N(x)+1)\,dx.
\ref{\thelaplaceref}
$$

\thm\classlimitthm{In addition to the assumptions above, let $\phi(s)$ 
satisfy conditions (1)--(3) and ($\alpha$)--($\gamma$) of Appendix
{\tauberapp}. Then:
\item{1.}  For every $f\in C(\ol\bR_\Omega)$,
$$
\lim_{E\to\infty}\tau_E(T(f))=\int_{ap}f(x)\,dx.
$$
\item{2.} For every $f,g\in C(\ol\bR_\Omega)$,
$$
\lim_{E\to\infty}\tau_E(T(f)T(g)-T(fg))=0.
$$
\item{3.} For every $a\in\cA_+$ the limit
$$
\lim_{E\to\infty}\tau_E(a)=:\tau(a)
$$ 
exists.
\item{4.} Let $\pi_\tau$ be the GNS representation of $\cA_+$ with
respect to the state $\tau$. Then, 
$$
\pi_\tau(\cA_+)\simeq C(\ol\bR_\Omega),
$$
and for every $a\in\cA_+$ we have 
$$
\pi_\tau(U(t)aU(-t))=
\alpha_t(\pi_\tau(a)).
$$}

\medskip\noindent{\it Proof.}
We first introduce some notation. As before, $u_\omega=T(e^{ix_\omega})$
denotes the unilateral shift. We set 
$$
u_\omega(n):=\cases{u_\omega^n, &if $n\geq 0$,\cr
(u_\omega^*)^{-n}, &if $n<0$.}
$$
It is easy to verify that
$$
U(t)u_\omega(n)U(-t)=e^{in\omega t}u_\omega(n).\ref{\uuref}
$$
If $f\in C(\ol\bR_\Omega)$ is a trigonometric polynomial,
$$
f\big(\prod e^{ix_\omega}\big)=\sum_{\{n_\omega\}}f_{n_\omega}\,
e^{in_\omega x_\omega},
$$
then the corresponding Toeplitz operator $T(f)$ is explicitly given by
$$
T(f)=\sum_{n_\omega}f_{n_\omega}\,u_\omega(n_\omega).
$$
It follows that for any $\eta$,
$$
(e(\eta),T(f)e(\eta))=f_{0}=\int_{ap}f(x)\,dx.
$$
This and the continuity of $T(f)$ in $f$ prove part 1 of \classlimitthm.

For later reference, notice also that as a consequence of \uuref,
$$
U(t)T(f)U(-t)=T(\alpha_t(f)).\ref{\alpharef}
$$
 
The structure of the standard Toeplitz algebra implies that operators 
$T(f)T(g)-T(fg)$ generate the commutator ideal $\cal{I}$ of $\cA_+$.
The quotient $\cA_+/\cal{I}$ is isomorphic to $C(\ol\bR_\Omega)$, and
the quotient map $\pi:\cA_+\to C(\ol\bR_\Omega)$ is called the symbol map.
We claim that
$$
\tau(a)=\lim_{E\to\infty}\tau_E(a)=\int_{ap}\pi(a)(x)\,dx,.
\ref{\taustateref}
$$
for all $a\in\cA_+$. In other words, the state $\tau$ is trivial on
the commutator ideal, and it coincides with the Lebesgue integral on
the abelian quotient. Parts 2 and 3 of \classlimitthm\ are
straightforward consequences of \taustateref . Formula \taustateref\
implies also that $\pi_\tau(\cA_+)$ is isomorphic the algebra
$C(\ol\bR_\Omega)$. The last statement of the theorem follows now
from \alpharef . 

To prove \taustateref , it is enough, in view of part 1 of \classlimitthm , 
to show that
$$
\lim_{E\to\infty}\tau_E(a)=0,\ \ {\rm if}\ \ a\in\cal{I}\ .
$$
The structure of the standard Toeplitz algebra $\fT$ implies that
$\cal{I}$ is generated by the operators of the form
$$
a=a_{\omega_1}\otimes a_{\omega_2}\otimes\ldots\otimes a_{\omega_N}
\otimes I\otimes I\ldots,\ref{\optyperef}
$$
where at least one of the operators $a_{\omega_1}\ldots a_{\omega_N}$,
say $a_{\omega_k}$, is compact. By a density argument, it is no loss
of generality to assume that $a_{\omega_k}$ is a finite rank operator
whose range is spanned by finitely many elements of the canonical basis.
Let $\Pi$ be the orthogonal projection onto this subspace. Then
$$
\tau_E(a)\leq {\Vert a\Vert \over N(E)}\sum_{\bN[\Omega_+]\ni\eta\leq E}
(e(\eta),\Pi e(\eta))\ .
$$
Since the spectrum of $H_+$ is the set $\{\sum n_\omega\omega :
n_\omega\geq 0\}$, we have to show that for any integer $M$,
$$
{\#\{n_\omega\geq 0\ ,\ n_{\omega_k}\leq M : \sum n_\omega\omega\leq E\}
\over\#\{n_\omega\geq 0 : \sum n_\omega\omega\leq E\}}
\longrightarrow 0,\ {\rm as}\ E\to 0\ .\ref{\combref}
$$
The numerator of the LHS of \combref\ is equal to
$N(E)-N(E-\omega_k(M+1))$ and so we have to show that
$$
{N(E)-N(E-\omega_k(M+1))\over N(E)}\longrightarrow 0,\ {\rm as}\ E\to 0.
$$ 
Formula {\thelaplaceref} and Corollary A.3 yield $N(E+1)=N(E)\big(1+o(1)
\big)$. Consequently, \combref\ follows, and the
theorem is proved.
$\square$

\medskip\noindent
{\bf Remark}. If $\Omega_+=\{\log p:\ p\ \hbox{prime number}\}$, then
it is easy to see that $N(E)\sim e^E$ and \combref\ is not true.
It would be interesting to determine the classical limit in this
case. Note also that the prime number theorem implies that
$\phi(s)$ is divergent for $s\leq 1$, and so the Tauberian
theorem of Appendix \tauberapp\  does not apply.

\medskip It follows now from the general results in [Z] that the
quantum Kronecker flow $(\cA_+,U(t))$ is quantum ergodic in the following
sense.
\thm\ergodicthm{Under the assumptions of \classlimitthm ,
for every $a\in\cA_+$,
$$
\lim_{M\to\infty}{1\over M}\int_0^M U(t)aU(-t)\,dt=\tau(a)I+A,
$$
where $A$ is in the weak closure of $\cA_+$, and
$$
\lim_{E\to\infty}\tau_E(A^*A)=0\ .
$$}

\noindent
In other words, the time average of a quantum observable is equal
to its spatial average plus a correction which vanishes in the
classical limit.

\section\bosesec{Almost periodic Bose field}
\medskip
In this section we define the free almost periodic quantized field. It
arises as the result of canonical quantization of the classical almost
periodic wave equation. Using Bohr's mean, we propose a canonical
formulation of the latter, and apply the standard quantization procedure.
The resulting quantum dynamical system is a ``double" of the Kronecker
dynamics studied in previous section.  
\newsubsec
We first introduce some notation. For a smooth function $H:\cA\cP(\Omega)
\to\bC$ we let $\nabla_f H$ denote the Frechet derivative of $H$ in the
direction of $f\in \cA\cP(\Omega)$. The functional derivative ${\delta
H(\phi)\over\delta\phi(x)}$ is, by definition, the distribution on
$\cA\cP(\Omega)$ such that
$$
\nabla_fH(\phi)=\int_{ap}{\delta H(\phi)\over\delta\phi(x)}f(x)\,dx\ .
$$
The canonical complex structure on $\cA\cP(\Omega)$ defines a symplectic 
structure on $\cA\cP(\Omega)$ for which positions are real functions
and momenta are purely imaginary functions. The respective coordinates
will be denoted by $\phi(x)$ and $\pi(x)$ so that ${\delta H(\phi,\pi)\over
\delta\phi(x)}$ is the functional derivative in the real direction,
and ${\delta H(\phi,\pi)\over\delta\pi(x)}$ is the functional derivative
in the imaginary direction. The symplectic space $\cA\cP(\Omega)$ is
the phase space for almost periodic field theory. 

Every smooth function $H:\cA\cP(\Omega)\to\bR$ defines a Hamiltonian
flow on $\cA\cP(\Omega)$ by
$$
\eqalign{
{d\phi(x,t)\over dt}&={\delta H\over\delta\pi(x)},\cr
{d\pi(x,t)\over dt}&=-{\delta H\over\delta\phi(x)}\ .\cr
}\ref{\hamilref}
$$
The Poisson bracket of two functions $F$ and $G$ on $\cA\cP(\Omega)$ is
defined by
$$
\{F,G\}=\int_{ap}\left({\delta F\over\delta\phi(x)}
{\delta G\over\delta\pi(x)}- {\delta F\over\delta\pi(x)}
{\delta G\over\delta\phi(x)}\right)\,dx
$$
and so the flow \hamilref\ can be written as
$$
{dF(\phi(t),\pi(t))\over dt}=\{F,H\}\ .\ref{\obserref}
$$
A straightforward calculation shows that
$$
\eqalign{
&\{\phi(x,t),\pi(y,t)\}=\delta_\Omega(x-y), \cr
&\{\phi(x,t),\phi(y,t)\}=0, \cr
&\{\pi(x,t),\pi(y,t)\}=0 \ .\cr
}\ref{\bracketref}$$

We will now formulate the almost periodic free field theory. The
dynamics is given by the wave equation,
$$
{\partial^2\phi\over \partial t^2}-
{\partial^2\phi\over \partial x^2}=0\ .\ref{\waveref}
$$
Equation \waveref\ can be written in the form \hamilref\ with the
Hamiltonian
$$
H(\phi,\pi)={1\over 2}\int_{ap}\big(\pi(x)^2+
(\partial_x\phi(x))^2\big)\,dx.
$$
For this Hamiltonian, equations \hamilref\ read
$$
{d\phi\over dt}=\pi,\ \ {d\pi\over dt}=\partial^2_x\phi,
\ref{\motionref}
$$
and lead to {\waveref}. The most general solution of \motionref\ can
be written in the following form
$$
\eqalign{
&\phi(x,t)=\sum_{\omega\in\Omega}(\phi_{1,\omega}e^{i\omega(x+t)}+
\phi_{2,\omega}e^{i\omega(x-t)}), \cr
&\pi(x,t)=\sum_{\omega\in\Omega}i\omega(\phi_{1,\omega}e^{i\omega(x+t)}-
\phi_{2,\omega}e^{i\omega(x-t)})\ .\cr
}\ref{\solutionref}
$$
The reality condition implies $\ol\phi_{1,\omega}=\phi_{1,-\omega}$
and $\ol\phi_{2,\omega}=\phi_{2,-\omega}$. In quantum field theory
it is convenient to parameterize the above solutions slightly
differently. We set
$$
a_\omega=\sqrt 2\,|\omega|^{1/2}\times\;
\cases{\phi_{1,-\omega}\  \ {\rm if}\ \ \omega>0,\cr
\phi_{2,-\omega}\  \ {\rm if}\ \ \omega<0\ .\cr
}$$
The new variables $a_\omega$ have the following Poisson brackets
$$\eqalign{
&\{a_\omega,\ol a_{\omega'}\}=\delta_{\omega,{\omega'}}, \cr
&\{a_\omega,a_{\omega'}\}=\{\ol a_\omega,\ol a_{\omega'}\}=0 \ .\cr
}\ref{\creabracketsref}$$
Equations \solutionref\ can then be recast in the following form:
$$
\eqalign{
&\phi(x,t)={1\over\sqrt 2}\sum_{\omega\in\Omega}|\omega|^{-1/2}
(\ol a_{\omega}e^{it|\omega|}+
a_{-\omega}e^{-it|\omega|})e^{i\omega x}, \cr
&\pi(x,t)={1\over\sqrt 2}\sum_{\omega\in\Omega}|\omega|^{1/2}
(\ol a_{\omega}e^{it|\omega|}-
a_{-\omega}e^{-it|\omega|})e^{i\omega x}\ .\cr
}\ref{\fieldsref}$$

\subsec
We shall now describe a quantization of the algebra of functions 
$\cA\cP(\Omega)$ and of the dynamics {\motionref}. We will follow
the procedure of canonical quantization which is adopted in quantum
field theory.

The standard rule of quantization consists in replacing classical
observables by operators and Poisson brackets by 
${1\over i}\times$commutators. For simplicity we set $\hbar=1$.
More precisely, quantization of the almost periodic wave equation
proceeds as follows. We find almost periodic, hermitian, operator
valued distributions $\phi(x,t)$ and $\pi(x,t)$ such that \motionref\
is satisfied. Furthermore, we require that (see \bracketref)
$$
\eqalign{
&[\phi(x,t),\pi(y,t)]=i\delta_\Omega(x-y), \cr
&[\phi(x,t),\phi(y,t)]=0, \cr
&[\pi(x,t),\pi(y,t)]=0.\cr
}\ref{\commutatorref}
$$
The quantum hamiltonian $H_b$ determines the time evolution
of the field operators given by the Heisenberg equations of motion,
$$
{d\phi\over dt}={1\over i}[\phi,H_b],\ 
{d\pi\over dt}={1\over i}[\pi,H_b].\ref{\hamildefref}
$$
We construct operators $a_\omega$ and $a_\omega^*$ 
as in \fieldsref , satisfying the commutation relations
$$\eqalign{
&[a_\omega,a_{\omega'}^*]=\delta_{\omega,{\omega'}}, \cr
&[a_\omega,a_{\omega'}]=[a_\omega^*,a_{\omega'}^*]=0,\cr
}\ref{\creacommref}
$$
and such that $a_\omega^*$ is the hermitian conjugate of $a_\omega$.
Operators $a_\omega$ and $a_\omega^*$ are called annihilation and
creation operators, respectively.

If one additionally assumes the existence of a cyclic vector $v_0$
such that $a_\omega v_0=0$ and some natural domain restrictions,
it is known that the algebra {\creacommref} has a unique representation
in terms of a Fock space which we will describe now.

Let  $\cH=l^2(\Omega)$ and consider the bosonic Fock space $\cF_b\cH$.
As before, the vector $v_0=(1,0,0,\ldots)\in \cF_b\cH$ is
called the vacuum. The Hilbert space $\cF_b\cH$ can be naturally
identified with $l^2(\bN[\Omega])$, where $\bN[\Omega]$ is the
set of nonnegative, integer, finite combinations of elements of $\Omega$.
Alternatively,
$$
\cF_b\cH\simeq\bigotimes_{\omega\in\Omega}l^2(\bN[\omega]),
\ref{\bosprodref}
$$
as in \kronprodref .

Let $e(\eta)$, $\eta\in\bN[\Omega]$ be the canonical
orthonormal basis in $l^2(\bN[\Omega])\simeq\cF_b\cH$.
The set $\bN[\Omega]$ is, in a natural way, a semigroup with respect
to addition. Writing
$$
\bN[\Omega]\ni\eta=\sum n_\omega\omega,\ \omega\in\Omega, 
$$
where almost all numbers $n_\omega$ are zero, we define
the creation operators $a_\omega^*$ by
$$
a_\omega^* e(\eta)=\sqrt{n_\omega+1}\,e(\eta+\omega).\num
$$
The field operators $\phi(x,t)$ and $\pi(x,t)$ are then defined by means
of formula \fieldsref .

The Hamiltonian of the free almost periodic quantum field
theory is given by the familiar expression
$$
H_b=\sum_{\omega\in\Omega}|\omega| a_\omega^* a_\omega=
{1\over 2}\int_{ap}:\big(\pi(x)^2+(\partial_x\phi(x))^2\big):\,dx\ ,
\num
$$
where $:\;\; :$ means Wick ordering. The canonical basis $\{e(\eta)\}$
is the basis of eigenvectors for $H_b$,
$$
H_be(\eta)=(\sum n_\omega|\omega|)e(\eta).
\ref{\boshamdefref}
$$

Let $\cD$ be the dense subspace of $\cF_b\cH$ consisting of
finite linear combinations of the basis elements $e(\eta)$.
It is an invariant domain for $a_\omega$ and $a_\omega^*$, and 
is a core for $H_b$.

\prop\quantizationprop{With the above definitions, the operator
valued distributions $\phi(x,t)$ and $\pi(x,t)$, and the
Hamilton operator $H_b$ satisfy equations \motionref ,
\commutatorref\ and \hamildefref\ on $\cD$.
}

\medskip\noindent{\it Proof.} The proof is a direct calculation
following essentially the similar argument in standard quantum field
theory, see e.g. [GJ].
$\square$

\medskip The quantum dynamics described in this section is very
closely related to the quantum Kronecker flow. Indeed, denoting
$\cH_-:=l^2(\Omega_-)$, we have a natural decomposition
$$
\cF_b\cH\simeq\cF_b\cH_-\otimes\cF_b\cH_+.
$$
With respect to this decomposition the Hamiltonian $H_b$ can be split
into the positive and negative frequency parts:
$$
H_b=\sum_{\omega\in\Omega}|\omega| a_\omega^* a_\omega=
\sum_{\omega\in\Omega_-}|\omega| a_\omega^* a_\omega+
\sum_{\omega\in\Omega_+}\omega a_\omega^* a_\omega=
H_-+H_+\ .
$$
Since $H_-$ is unitarily equivalent to $H_+$, the almost periodic free
field theory is a double of the quantum Kronecker flow.

Consider the family of operators $U_\omega(t)$ and $V_\omega(s)$,
$s,t\in\bR$, $\omega\in\Omega$, defined as
$$
U_\omega(t)=e^{it(a_\omega+a_\omega^*)}\ ,\ 
V_\omega(s)=e^{s(a_\omega-a_\omega^*)}\ .
$$
The $\bC^*$-algebra generated by $U_\omega(t)$ and $V_\omega(s)$ is
an example of a CCR algebra [BR]. It is a nonseparable $C^*$-algebra
which is usually studied in quantum field theory. It is, however, not well
suited for our purposes and, following [BC], we define
the bosonic algebra of observables $\cA_b$ to be the $C^*$-algebra
generated by the canonical unilateral shifts in each factor of 
\bosprodref . The hamiltonian $H_b$ defines a dynamics $\sigma^b_t$
on $\cA_b$ by
$$
\sigma^b_t(A)=e^{itH_b}Ae^{-itH_b}\ ,\ A\in\cA_b\ .
$$
We have the following
analog of \toeplitzprop .

\prop\bosKMSprop{
\item{1.} The $\bC^*$-algebra $\cA_b$ is an infinite tensor
product of standard Toeplitz $\bC^*$-algebras $\fT_\omega$:
$$
\cA_b=\bigotimes_{\omega\in\Omega}\fT_\omega,
$$
where $\fT_\omega$ is generated by the canonical unilateral shift in
$l^2(\bN[\omega])$.
\item{2.} For every $\beta>0$, there exists a unique
KMS$_\beta$ state on $(\cA_b,\sigma^b_t)$.
}

\medskip\noindent{\it Proof.} This follows from Propositions 7 and 8 of
[BC].
$\square$

\section\fermsec{Almost Periodic Fermions}
\medskip
In this section we will define the almost periodic fermionic quantum free
field.
\newsubsec
Let, as before, $\cH=l^2(\Omega)$, and consider the fermionic Fock space
$\cF_f\cH$. The Hilbert space $\cF_f\cH$ is defined as
$$
\cF_f\cH=\bigoplus_{n=0}^{\infty}\bigwedge\nolimits^n\cH,
$$
where $\bigwedge^n\cH$ is the $n^{\rm th}$ exterior power of $\cH$
with $\bigwedge^0\cH=\bC$. The vector $v_0=(1,0,0,\ldots)\in \cF_f\cH$ is
called the vacuum. The Hilbert space $\cF_f\cH$ can be naturally
identified with $l^2(\bZ_2[\Omega])$, where $\bZ_2$ is
the group $\{0,1\}$ with addition modulo $2$.

Let $f(\eta)$, $\eta\in\bZ_2[\Omega]$ be the canonical
orthonormal basis in $l^2(\bZ_2[\Omega])\simeq\cF_f\cH$.
The set $\bZ_2[\Omega]$ has a natural group structure with
respect to addition modulo $2$. Writing
$$
\bZ_2[\Omega]\ni\eta=\sum_{\omega\in\Omega} n_\omega\omega, 
$$
where almost all numbers $n_\omega$ are zero, we define
the creation operators $b_\omega^*$  and the annihilation operators 
$b_\omega$ by
$$\eqalign{
&b_\omega^* f(\eta)=\sqrt{(n_\omega+1)\,{\rm mod}\,2}\ f(\eta+\omega),\cr
&b_\omega f(\eta)=\sqrt{n_\omega}\ f(\eta-\omega).\cr
}
$$
It easy to verify the following anticommutation relations
$$\eqalign{
&[b_\omega,b_{\omega'}^*]_+=\delta_{\omega,{\omega'}},\cr
&[b_\omega,b_{\omega'}]_+=[b_\omega^*,b_{\omega'}^*]_+=0,\cr
}\ref{\anticommref}
$$
where $[x,y]_+:=xy+yx$ is the anticommutator. Unlike in
the bosonic case, the operators $b_\omega$ are bounded.
Let $\cA_f$ be the $\bC^*$-algebra generated the fermionic 
creation and annihilation operators. This algebra is called in the
literature the CAR algebra [BR].

Fermionic field operators $\psi_1(x)$ and $\psi_2(x)$ at time 0
are then defined by
$$\eqalign{
&\psi_1(x)=\sum_{\omega\in\Omega}\left(\Theta(\omega)b_\omega^*+
\Theta(-\omega)b_{-\omega}\right)e^{-i\omega x},\cr
&\psi_2(x)=i\sum_{\omega\in\Omega}\left(\Theta(-\omega)b_\omega^*+
\Theta(\omega)b_{-\omega}\right)e^{-i\omega x},\cr
}\ref{\fermionopref}
$$
where $\Theta$ is the Heaviside function. One can directly
verify the following anticommutation relations
$$
[\psi_i(x),\psi_j(y)]_+=2\delta_{ij}\delta_\Omega(x-y)\ .
$$
The fermionic almost periodic free hamiltonian $H_f$ is then 
given by
$$
H_f=\sum_{\omega\in\Omega}|\omega| b_\omega^* b_\omega\ ,
$$
so that
$$
H_ff(\eta)=(\sum_{\omega\in\Omega} n_\omega|\omega|)f(\eta).
$$
It defines a dynamics $\sigma^f_t$ on $\cA_f$ by
$$
\sigma^f_t(A)=e^{itH_f}Ae^{-itH_f}\ ,\ A\in\cA_f.
$$

\subsec
The supersymmetric almost periodic quantum free field theory
is defined as the tensor product of bosonic and fermionic 
field theories. This means that the Hilbert space of that
theory is the tensor product $\cF_b\cH\otimes\cF_f\cH$
with the natural $\bZ_2$ grading $\Gamma=I\otimes(-1)^F$,
where
$$
Ff(\eta)=(\sum_{\omega\in\Omega} n_\omega)f(\eta).
$$
The relevant $\bC^*$-algebra $\cA$ is then the tensor product
$\cA=\cA_b\otimes\cA_f$. The supersymmetric hamiltonian $H$
is defined by
$$
H=H_b\otimes I+I\otimes H_f,
$$
and the corresponding dynamics on $\cA$ is denoted by $\sigma_t$.
The new feature of the supersymmetric theory is the existence 
of a supercharge, namely a self-adjoint operator $Q$ which is odd
under the $\bZ_2$-grading, and has the property that $Q^2=H$.
The operator $Q$ can be defined in the following way:
$$\eqalign{
Q&={1\over\sqrt 2}\int_{ap}\psi_1(x)\left(
(\pi(x)-\partial_x\phi(x))+\psi_2(x)(\pi(x)+\partial_x\phi(x))
\right)\,dx\cr
&=\sum_{\omega\in\Omega}\sqrt{|\omega|}(a_\omega^*b_\omega+
a_\omega b_\omega^*)\ .\cr
}\ref{\Qdefref}
$$
The system $(\cA,\Gamma,\sigma_t,Q)$ is an example of a
quantum algebra to be discussed in the next section.

\section\skmssec{Super KMS States}
\medskip
In this section we construct and prove the uniqueness of super-KMS
functionals for the free supersymmetric almost periodic quantum field
theory. Super-KMS functionals are $\bZ_2$-graded counterparts of
KMS states and play an important role in index theory.
\newsubsec
We will recall the definitions of quantum algebras and the super-KMS
states on quantum algebras [K], [JLW]. 

\defin\quantalgdefine{A {\it quantum algebra} is a quadruple
$(\cA,\Gamma,\sigma_t,d)$ satisfying conditions 1-4 below.
\item{1.} {$\cA$ is a $\bC^*$-algebra.}
\item{2.} {$\Gamma$ is a $\bZ_2$ grading on $\cA$ i.e. a $*$-automorphism
of $\cA$ such that $\Gamma^2=I$. For $a\in\cA$ we denote 
$a^\Gamma:=\Gamma(a)$.}
\item{3.} {$\sigma_t:\cA\to\cA$ is a continuous, one-parameter group of
even, bounded automorphisms of $\cA$. $\sigma_t$ do not have to be
$*$-automorphisms.}
\item{4.} {Let $\cA_\alpha$ be the subalgebra of $\cA$ such that
for every $a\in\cA_\alpha$ the function $t\to\sigma_t(a)$
extends to an entire $\cA$-valued function. It is known that
$\cA_\alpha$ is norm dense. On $\cA_\alpha$ we set
$$
D:=-i{d\sigma_t\over dt}\big|_{t=0}
.$$
$d$ is a superderivation on $\cA_\alpha$ i.e.
$$
d^\Gamma=-d\ ,\ d(ab)=da\,b+a^\Gamma\,db,
$$
such that $d^2=D$.}
\item{}}

In the theory of the previous section set $da:=[Q,a]_s$,
and $Da:=[H,a]_s$, where $[a,b]_s$ is the supercommutator, i.e.
$[a,b]_s=[a,b]$, if at least one of the operators $a, b$ is even,
and $[a,b]_s=[a,b]_+$, if both are odd.
Then $(\cA,\Gamma,\sigma_t,d)$ is a quantum algebra. In the
following, this quantum algebra will be referred to as
the {\it almost periodic quantum algebra}.

\defin\skmsdefine{Let $(\cA,\Gamma,\sigma_t,d)$ be a quantum algebra.
A continuous linear functional $\mu:\cA\to\bC$ is called
a super-KMS$_\beta$ functional if for $a,b\in\cA_\alpha$,
\item{1.} $\mu(da)=0$,
\item{2.} $\mu(ab)=\mu(b^\Gamma\sigma_{i\beta}(a))$.
\smallskip\noindent
If, for a $\bZ_2$-graded $\bC^*$ dynamical
system $(\cA,\Gamma,\sigma_t)$, a linear continuous functional
$\mu$ satisfies only the condition 2 above, then it is called
a {\it pre super-KMS$_\beta$ functional}.}
\medskip
Unlike for KMS states, no positivity assumptions
are or can be made for super-KMS$_\beta$ functionals. 
It follows from the definition that a super-KMS$_\beta$
functional $\mu:\cA\to\bC$ satisfies
$$\eqalign{
&\mu(\sigma_t(a))=\mu(a),\cr
&\mu^\Gamma=\mu,\cr
&\mu(a\,db)=\mu(da\,b^\Gamma).\cr
}$$
In our example of the almost periodic quantum algebra, 
assuming additionally that $\tr(e^{-\beta H})<\infty$, set
$$
\mu_\beta(a):=\str(ae^{-\beta H})\ ,\ref{\canskmsdef}
$$
where $\str$ is the supertrace, i.e. $\str(a)=\tr(\Gamma a)$. It is
then easy to verify that $\mu_\beta$ is a super-KMS$_\beta$ functional.

\subsec
The remainder of this section is devoted to the proof of the uniqueness
of the super-KMS functional for the almost periodic quantum algebra. We
start with two propositions of independent interest.

\prop\finitedimprop{Let $V$ be a finite dimensional Hilbert space,
$\Gamma$ a $\bZ_2$-grading on $V$, $Q$ an odd self-adjoint operator
on $V$, $H:=Q^2$, and $\cA:=\cL(V)$ the algebra of linear operators on $V$.
For $a\in\cA$, we set $da:=[Q,a]_s$ and $\sigma_t(a):=e^{itH}ae^{-itH}$.
Then, for every $\beta>0$, there is a unique, up to a 
multiplicative constant, 
pre super-KMS$_\beta$ functional $\mu_\beta$ on $(\cA,\Gamma,\sigma_t)$ 
given by
$$
\mu_\beta(a)=\str(ae^{\beta H}).
$$
Moreover, $\mu_\beta$ is automatically a super-KMS$_\beta$ functional
on $(\cA,\Gamma,\sigma_t,d)$.
}
\medskip\noindent{\it Proof.}
Let $\mu_\beta$ be any pre super-KMS$_\beta$ functional on 
$(\cA,\Gamma,\sigma_t,d)$. 
Consider $\tilde\mu_\beta:=\mu_\beta(\Gamma ae^{-\beta H})$.
Using condition 2 of \skmsdefine , one easily verifies that
$\tilde\mu_\beta(ab)=\tilde\mu_\beta(ba)$, and so
$\tilde\mu_\beta$ is proportional to the trace and the claim follows. 
$\square$

\prop\tensprodprop{Let $(\cA^i,\Gamma^i,\sigma_t^i)$, $i=1,2$
be two $\bZ_2$-graded $\bC^*$-dynamical systems which have unique, up
to a multiplicative constant, pre super-KMS$_\beta$ functionals
$\mu_\beta^i$. Then $\mu_\beta^1 \otimes\mu_\beta^2$ is a
unique, up to a multiplicative constant, pre super-KMS$_\beta$ 
functional on the tensor product 
$(\cA^1\otimes\cA^2,\Gamma^1\otimes\Gamma^2,\sigma_t^1\otimes\sigma_t^2)$.
}
\medskip\noindent{\it Proof.}
Let $\mu_\beta$ be any pre super-KMS$_\beta$ functional on the 
tensor product.
The statement follows easily from the fact that for any $b\in\cA^2$
the following functional on $\cA^1$:
$$
\mu_{\beta,b}(a):=\mu_\beta(a\otimes b)
$$
is a pre super-KMS$_\beta$ functional.
$\square$
\medskip

The following theorem can now be easily deduced from \finitedimprop ,
\tensprodprop , and Proposition 8 of [BC].

\thm\sKMSthm{For every $\beta>0$, there exists a unique, up to a 
multiplicative constant,
super-KMS$_\beta$ functional on the almost periodic quantum
algebra $(\cA,\Gamma,\sigma_t,d)$.
}
\medskip\noindent{\it Proof.}
We are going to prove that there is a unique pre super-KMS$_\beta$ 
functional on the $\bZ_2$-graded $\bC^*$-dynamical system
$(\cA,\Gamma,\sigma_t)$. It will follow from the construction
that that functional is, in fact, a super-KMS$_\beta$ 
functional.

It follows from \bosKMSprop\ and the structure theorems for $\cA_f$,
see [BR], that the $\bC^*$-algebra $\cA$ is isomorphic with
the following infinite tensor product:
$$
\cA=\bigotimes_{\omega\in\Omega}\fT_\omega\otimes\fA_\omega,
$$
where $\fT_\omega$ is the Toeplitz algebra and
$\fA_\omega$ is generated by the fermionic creation
and annihilation operators $b_\omega^*,b_\omega$, and is isomorphic 
with $M_2(\bC)$, the algebra of $2\times 2$ matrices. Additionally, 
both the grading $\Gamma$ and the dynamics $\sigma_t$ factor
with respect to the above decomposition,
$$
\Gamma=\bigotimes_{\omega\in\Omega}\Gamma_\omega,\quad 
\sigma_t=\bigotimes_{\omega\in\Omega}\sigma_{t,\omega}.
$$
It is easy to verify that $\Gamma_\omega$ is trivial on $\fT_\omega$,
so that $\Gamma_\omega=I\otimes\Gamma_\omega^f$. The generator 
of $\sigma_{t,\omega}$ is the supersymmetric harmonic
oscillator hamiltonian $H_\omega=|\omega|(a_\omega^*a_\omega
+b_\omega^*b_\omega)$, and so we have a further decomposition: 
$\sigma_{t,\omega}=\sigma_{t,\omega}^b\otimes\sigma_{t,\omega}^f$.
The system $(\fA_\omega,\Gamma^f,\sigma_{t,\omega}^f)$ is finite 
dimensional, and thus by \finitedimprop\ it has a unique pre
super-KMS$_\beta$ functional. The uniqueness of a pre super-KMS$_\beta$ 
functional on $(\fT_\omega,\Gamma^b=I,\sigma_{t,\omega}^b)$
follows from Proposition 8 of [BC] since the proof of that 
proposition does not require any positivity assumptions on the functional.
Moreover, any pre super-KMS$_\beta$ functional on
$\fT_\omega\otimes\fA_\omega\ $ is proportional to
$$
a\to\str(ae^{-\beta H_\omega}),
$$
and consequently is a super-KMS$_\beta$ functional.
The theorem now follows from \tensprodprop .
$\square$

\appendix\tauberapp{An Ingham type tauberian theorem}
\medskip
In this appendix we prove a technical result used in Section {\kronsec}
to establish the quantum ergodicity of the quantized Kronecker
dynamics. This result is a variant of Ingham's tauberian theorem [I],
see also [P], and differs from the original theorem in some of the
hypotheses.

Let $N(x)$ be a nondecreasing function of bounded variation
satisfying the following assumptions:
\item{(1)} {$N(x)=0$, for all $x\leq 0$;}
\item{(2)} {for all $\sigma>0$, $\int_0^\infty e^{-\sigma x}dN(x)
<\infty$;}
\item{(3)} {for all $s=\sigma+it$, $\sigma>0$, $t\in\bR$, the
function $\phi(s)$ defined by 
$$
e^{\phi(s)}=\int_0^\infty e^{-sx}dN(x) \ref{\lapltransf}
$$
is holomorphic.}

\noindent
The function $\phi(s)$ will play a fundamental role in the following
analysis. Ingham's original theorem requires detailed knowledge
of the asymptotic of $\phi(s)$ as $s$ approaches $0$ within an angle.
Such an asymptotic is usually difficult to obtain. Somewhat different
assumptions on $\phi(s)$ lead to a result which is well tailored for
our purposes. Specifically, we require that:
\item{($\alpha$)} {
$$
-\sigma\phi^\prime(\sigma)\nearrow\infty,\quad\hbox{and }
\sigma^2\phi^{\prime\prime}(\sigma)\nearrow\infty, \hbox{ as }\sigma
\searrow 0;\ref{\alpharef}
$$}
\item{($\beta$)} {
$$
\big|{{\sigma\phi^{\prime\prime\prime}(\sigma)}\over
{\phi^{\prime\prime}(\sigma)}}\big|=O(1), \hbox{ as }
\sigma\searrow 0;\ref{\betaref}
$$}
\item{($\gamma$)} {for any $\Delta>0$, there is $\sigma_0>0$ such
that the triangle
$$
T(\Delta,\sigma_0)=\{\sigma+it:\;0<\sigma<\sigma_0,\;|t|<\Delta\sigma\}
$$
does not contain nonreal roots of $\im\phi^{\prime}(s)$.}

We can now formulate the main result of this appendix..
\thm\tauberianthm{Under the above assumptions (1 - 3) and
($\alpha$ - $\gamma$),
$$
N(E)=\big(2\pi\sigma_E^2\phi^{\prime\prime}(\sigma_E)\big)^{-1/2}
e^{\sigma_EE+\phi(\sigma_E)}\big(1+o(1)\big), \hbox{ as }
E\rightarrow\infty,\ref{\asymptotic}
$$
where $\sigma_E$ is the unique solution to the equation
$$
\phi^\prime(\sigma)+E=0.\ref{\crpteq}
$$}
\proof
The existence and uniqueness of the solution of \crpteq\ follows from
assumption ($\alpha$).
Integrating by parts in the right hand side of {\lapltransf} we
obtain the identity
$$
{{e^{\phi(\sigma+it)}}\over{\sigma+it}}=\int_0^{\infty}e^{-(\sigma
+it)x}N(x)dx.\ref{\lapltransftwo}
$$
Let $g$ be an integrable function. Multiplying {\lapltransftwo}
by $e^{E(\sigma+it)}g(t)$ and integrating over $t$ we obtain,
after a change of the order of integration,
$$
\int_{-\infty}^\infty{{e^{\psi_E(\sigma+it)}}\over{\sigma+it}}
g(t)dt=\sqrt{2\pi}\int_{-\infty}^\infty e^{\sigma(E-x)}
\widehat g(x-E)N(x)dx,\ref{\laplintegrated}
$$
where we have set
$$
\psi_E(s)=\phi(s)+Es.\ref{\psidef}
$$
Shifting the integration variable in the right hand side of
{\laplintegrated} we rewrite {\laplintegrated} as the following
basic identity
$$
{1\over{2\pi}}\int_{-\infty}^\infty{{e^{\psi_E(\sigma+it)}}\over{\sigma+it}}
g(t)dt={1\over{\sqrt{2\pi}}}\int_{-\infty}^\infty e^{-\sigma x}
\widehat g(x)N(x+E)dx.\ref{\basicidentity}
$$
We take the function $g$ to be of the form $g(t)=f(t/T)$, where
$f$ is continuous in the interval $[-1,1]$ and zero outside it, $f(0)=1$,
and where $T>0$ is a number which will be chosen shortly. We let
$L(\sigma)$ denote the left hand side of {\basicidentity}, i.e.
$$
L(\sigma)={1\over{2\pi}}\int_{-T}^T{{e^{\psi_E(\sigma+it)}}\over{\sigma+it}}
f(t/T)dt.
$$
For $0<\delta<T$, we decompose $L(\sigma)$ into two parts
$$
\eqalign{
L(\sigma)&={1\over{2\pi}}\int_{|t|\leq\delta}{{e^{\psi_E(\sigma+it)}}\over
{\sigma+it}}f(t/T)dt+{1\over{2\pi}}\int_{\delta<|t|<T}
{{e^{\psi_E(\sigma+it)}}\over{\sigma+it}}f(t/T)dt\cr
&\equiv L^{(1)}(\sigma)+L^{(2)}(\sigma),\cr}
$$
and analyze them separately.

So far the considerations have been quite general, and we will now
start making specific choices. Pick any $\Delta>0$ (which
we will eventually want to make arbitrarily large), and choose
$\sigma_0>0$ such that the triangle $T(\Delta,\sigma_0)$ defined
in assumption ($\gamma$) does not contain nonreal roots of
Im$\phi^\prime(s)$. Take $E$ sufficiently large so that $\sigma_E
<\sigma_0$. To simplify the notation, $\sigma_E$ will be denoted by
$\sigma$ in the throughout the remainder of this proof. Furthermore, 
take $E$ large enough so that
$$
{\Delta\over\sqrt{\sigma^2\phi^{\prime\prime}(\sigma)}}\leq 1 ,
\ref{\Deltacond1}
$$
which is possible by assumption ($\alpha$). Set $T=\sigma\Delta$.
The choice of $\delta$ will be made shortly.

To analyze $L^{(1)}(\sigma)$ we expand $\psi_E(s)$
around $s=\sigma$ (in the following the subscript $E$ in $\psi_E$
will be suppressed):
$$
\psi(\sigma+it)=\psi(\sigma)-1/2\;\phi^{\prime\prime}(\sigma)t^2
-1/6\;\phi^{\prime\prime\prime}(\theta)it^3,
$$
for a $\theta$ belonging to the line segment which joins $\sigma-i\delta$
and $\sigma+i\delta$. By assumption ($\beta$) we have,
$$
|\phi^{\prime\prime\prime}(\theta)t|\leq\big|
{\sigma\phi^{\prime\prime\prime}(\sigma)\over\phi^{\prime\prime}(\sigma)}
\big|\,{\delta\over\sigma}\,\phi^{\prime\prime}(\sigma)
=o(1)\phi^{\prime\prime}(\sigma),
$$
if $\delta/\sigma=o(1)$, as $E\to\infty$. We now make the following
choice of $\delta$:
$$
\delta=\left({\sigma^2\over\phi^{\prime\prime}(\sigma)}\right)^{1/4}.
\ref{\deltachoice}
$$
Then, by assumption ($\alpha$), $\delta/\sigma=(\sigma^2\phi^{\prime\prime}
(\sigma))^{-1/4}=o(1)$, and consequently
$$
\psi(\sigma+it)=\psi(\sigma)-1/2\;\phi^{\prime\prime}(\sigma)(1+o(1))t^2.
$$
Therefore,
$$
\eqalign{
L^{(1)}(\sigma)&={1\over{2\pi}}\int_{-\delta}^\delta{1\over\sigma}f(0)
e^{\psi(\sigma)}e^{-1/2\;\phi^{\prime\prime}(\sigma)(1+o(1))t^2}\cr
&={e^{\psi(s)}\over{2\pi\sqrt{\sigma^2\phi^{\prime\prime}(\sigma)}}}
\int_{-\delta(\phi^{\prime\prime}(\sigma))^{1/2}}^
{\delta(\phi^{\prime\prime}(\sigma))^{1/2}}e^{-1/2\;(1+o(1))t^2}
\,dt\,\;(1+o(1)).\cr}
$$
But $\delta(\phi^{\prime\prime}(\sigma))^{1/2}=
(\sigma^2\phi^{\prime\prime}(\sigma))^{1/4}\to\infty$, as $E\to\infty$,
and so the gaussian integral above becomes an integral
over entire $\bR$ in this limit. As a result,
$$
L^{(1)}(\sigma)=\big(2\pi\sigma^2\phi^{\prime\prime}
(\sigma)\big)^{-1/2}e^{\psi(\sigma)}\big(1+o(1)\big).
$$

We now turn to the analysis of $L^{(2)}(\sigma)$. We wish to show
that this term is much smaller than the previous one. Indeed,
$$\eqalign{
\big|L^{(2)}(\sigma)\big|&\leq T\,\sup |f(t/T)|\,{1\over\sigma}
\,\sup_{\delta<|t|<T} \big|e^{\psi(\sigma+it)}\big|\cr
&=O(1)\,\Delta\sup_{\delta<|t|<T}e^{\re\psi(\sigma+it)}.\cr
}$$
Assumption ($\gamma$) implies that the above supremum is attained
at $|t|=\delta$. To see this, we consider the function
$$
t\rightarrow\re\psi(\sigma+it)=\re\phi(\sigma+it)+E\sigma.
$$
The critical points of this function satisfy
$$
0={d\over dt}\re(\phi(\sigma+it)+E\sigma)=\im\phi^\prime(\sigma+it).
$$
Hence there are no critical points in the interval $\delta<|t|<T=
\sigma\Delta$, and the function attains its maximum value at an
endpoint. Consequently, using a Taylor expansion as in the analysis of
$L^{(1)}(\sigma)$,
$$
\big|L^{(2)}(\sigma)\big|\leq O(1)\,\Delta\,e^{\re\psi(\sigma\pm i\delta)}
=O(1)\,\Delta\,e^{\psi(\sigma)-1/2\phi^{\prime\prime}(\sigma)
\delta^2(1+o(1))}.
$$
It is easy to see that, with our choices of
$\delta$ and $\Delta$, we have
$$
\Delta e^{-1/2\phi^{\prime\prime}(\sigma)\delta^2}\ll\big(\sigma^2
\phi^{\prime\prime}(\sigma)\big)^{-1/2}
$$
and so $L^{(2)}(\sigma)=o(1)\,L^{(1)}(\sigma)$. This concludes the
analysis of the left hand side of {\basicidentity}.

The asymptotic behavior of $L(\sigma)$ turns out to be independent of
the choice of function $f$. In the following we study the right hand
side of {\basicidentity} which will be denoted by $R(\sigma)$. We shall
make suitable choices of $f$ in order to get bounds on $N(E)$ from
above and from below.
\lemma\shiftlemma{Define
$$
\tilde N(E):=\big(2\pi\sigma^2\phi^{\prime\prime}(\sigma)\big)^{-1/2}
e^{\sigma E+\phi(\sigma)}.\ref{\ntildedefref}
$$
Then
$$
\tilde N(E+O(1)/\sigma)=\tilde N(E)\big(1+o(1)\big).
$$
}
\proof
Let $\sigma_1$ be the unique solution of the equation
$-\phi^\prime(\sigma_1)=E+O(1)/\sigma$. Taylor expanding $\phi^\prime$
around $\sigma$ yields
$$
\sigma_1=\sigma+{O(1)\over \sigma\phi^{\prime\prime}(\sigma)},
$$
since $1/\sigma\phi^{\prime\prime}(\sigma)\ll\sigma$.
It then follows readily that $\sigma_1=\sigma\big(1+o(1)\big)$.
In a similar fashion, we conclude that
$\phi^{\prime\prime}(\sigma_1)=\phi^{\prime\prime}(\sigma)
\big(1+o(1)\big)$, and $\psi(\sigma_1)=\psi(\sigma)+o(1)$. Inserting
these expressions into the definition of $\tilde N(E+O(1)/\sigma)$
completes the proof.
$\square$\medskip
\noindent{\it Choice 1}. Set
$$
f(t)=\cases{1-|t|,&if $|t|\leq 1$,\cr
0,&otherwise.\cr
}$$
The Fourier transform of $f$ is
$$
\widehat f(x)= {1\over \sqrt{2\pi}}\left({\sin(x/2)\over x/2}\right)^2,
$$
and thus the right hand side of \basicidentity\ is
$$
R(\sigma)={T\over 2\pi}\int_\bR e^{-x\sigma}
\left({\sin(Tx/2)\over Tx/2}\right)^2\,N(E+x)\,dx.\ref{\rchoice1}
$$
The integrand of \rchoice1\ is positive and so, for any $\Lambda$,
$$
\eqalign{R(\sigma)&\ge{T\over 2\pi}\int_{-\Lambda}^\Lambda e^{-x\sigma}
\left({\sin(Tx/2)\over Tx/2}\right)^2\,N(E+x)\,dx\cr
&\ge N(E-\Lambda)\,e^{-\sigma\Lambda}\,{1\over 2\pi}
\int_{-T\Lambda}^{T\Lambda}\left({\sin(x/2)\over x/2}\right)^2\,dx,\cr
}$$
where we have used the monotonicity of $N(x)$. Now take $\Lambda=1/(\sigma
\sqrt\Delta)$. With this choice, $\sigma\Lambda=1/\sqrt\Delta$, and 
the exponential term in the above formula tends to $1$, as $\Delta
\rightarrow\infty$. Similarly, $T\Lambda=\sqrt\Delta$ and the integral
over $(-T\Lambda,T\Lambda)$ tends to the integral (equal to $1$) over all
of $\bR$. This yields the inequality
$$
R(\sigma)\ge N(E-\Lambda)\big(1+o(1)\big).
$$
Replacing $E$ by $E-\Lambda$ and using \shiftlemma , we conclude that
$$
N(E)\leq\tilde N(E)\big(1+o(1)\big).
$$\medskip
\noindent{\it Choice 2}. Set
$$
f_1(t)=\cases{{1\over 2i\mu}\left(e^{i\mu|t|}-e^{-i\mu|t|}\right)
,&if $|t|\leq 1$\cr
0,&otherwise,\cr
}$$
where $\mu=2k\pi$, $0<k\in\bZ$. Let
$$
f_2(t)=\cases{1-|t|,&if $|t|\leq 1$,\cr
0,&otherwise,\cr
}$$
and take $f=f_1+f_2$. The Fourier transform of $f$ is
$$
\widehat f(x)= {1\over \sqrt{2\pi}}\left({\sin(x/2)\over x/2}\right)^2\,
{\mu^2\over \mu^2-x^2},
$$
and so $\widehat f(x)<0$, for $|x|>\mu$. Consequently,
$$
\eqalign{R(\sigma)&\leq{T\over 2\pi}\int_{|Tx|\leq\mu} e^{-x\sigma}
\left({\sin(Tx/2)\over Tx/2}\right)^2\,{\mu^2\over \mu^2-x^2}\,N(E+x)\,dx\cr
&\leq N(E+\mu/T)\,e^{\mu\sigma/T}\,{1\over 2\pi}\int_{-\mu}^\mu 
\left({\sin(x/2)\over x/2}\right)^2\,{\mu^2\over \mu^2-x^2}\,dx,\cr
}
$$
by monotonicity. Now take $k$ to be the integer part of $[\sqrt\Delta]$.
With this choice, the integral above tends to $1$, as $\Delta\rightarrow
\infty$. Also, $\mu\sigma/T\sim\Delta^{-1/2}\rightarrow 0$ and the
exponential term tends to $1$. Since $\mu/T\sim1/\sigma\sqrt\Delta$,
we can use {\shiftlemma} to replace $E$ by $E+\mu/T$. This yields
$$
N(E)\ge\tilde N(E)\big(1+o(1)\big),
$$
and concludes the proof of the theorem.
$\square$
\medskip
\cor\maincor{With the above assumptions we have
$$
N\big(E+O(1)\big)=N(E)\big(1+o(1)\big).
$$
}
\proof
This follows directly from {\tauberianthm} and {\shiftlemma}.
$\square$

\appendix\examplapp{Some examples of Kronecker systems}
\medskip
The theorem below provides a source of examples of Kronecker systems
satisfying the assumptions of {\tauberianthm}.
\thm\examplesthm{Let $\omega_n=An^\alpha(1+\mu_n)$, where
$A>0$ and $\alpha\geq 1$ are constant, and where $\mu_n=o(1)$,
as $n\rightarrow\infty$. Then $\phi(s)$ satisfies the assumptions
of {\tauberianthm}.}
\proof Assumptions (1) -- (3) of Appendix {\tauberapp} are clearly
satisfied, and so it is sufficient to verify assumptions ($\alpha$)
-- ($\gamma$).

Assumption ($\alpha$) is a consequence of the following equalities:
$$
\eqalign{
-\sigma\phi^\prime(\sigma)&=\sum_{n=1}^\infty f(\sigma\omega_n),\cr
\sigma^2\phi^{\prime\prime}(\sigma)&=\sum_{n=1}^\infty f(\sigma\omega_n),
\cr}
$$
where the functions $f$ and $g$ are given by
$$
f(x)={{xe^{-x}}\over{1-e^{-x}}},\qquad
g(x)={{x^2e^{-x}}\over{(1-e^{-x})^2}}.
$$
Since both $f(x)$ and $g(x)$ increase monotonically $1$ as $x\searrow 0$,
the claim follows.

To prove ($\beta$), we note that
$$
-\sigma^3\phi^{\prime\prime\prime}(\sigma)\leq\theta(\sigma),
$$
and
$$
\sigma^2\phi^{\prime\prime}(\sigma)\geq C_\epsilon
\theta((1-\epsilon)\sigma),\hbox{ for some }0<\epsilon<1,
$$
where $\theta(\sigma)$ is defined in {\thetasumref}. Since
$\omega_n=An^\alpha(1+\mu_n)$ implies that $\theta(\sigma)=C
\sigma^{-1/\alpha}(1+o(1))$, as $\sigma\rightarrow 0$, we conclude that
$$
\big|{{\sigma\phi^{\prime\prime\prime}(\sigma)}\over
{\phi^{\prime\prime}(\sigma)}}\big|\leq O(1)\; {{\theta(\sigma)}\over
{\theta((1-\epsilon)\sigma)}}=O(1),
$$
as $\sigma\rightarrow 0$.

Finally, assumption ($\gamma$) is verified in the following lemma.
$\square$
\lemma\nozeroeslemma{Under the assumptions of {\examplesthm},  
$$
\im\phi^\prime(\sigma+ix\sigma)=C_{A,\alpha}\sigma^{-\beta}
x\big(h(x)+o(1)\big),\ref{\assymptref}
$$
for $x\in\bR$, as $\sigma\rightarrow 0$, uniformly in $|x|\leq\Delta$.
Here $\beta=1+\alpha^{-1}$, and $h(x)$ is a function such that $h(x)\neq 0$,
for all $x\in\bR$.}
\proof Explicitly,
$$
\im\phi^\prime(\sigma+ix\sigma)=\sum_{n\geq 1}
{{\lambda_n\,e^{-\sigma\lambda_n}\sin(x\sigma\lambda_n)}\over
{1+e^{-2\sigma\lambda_n}-2e^{-\sigma\lambda_n}\cos(x\sigma\lambda_n)}}.
\ref{\imphiprimeref}
$$
We will analyze this expression in two steps.

\noindent {\it Step 1}. Assume first that $\lambda_n=An^\alpha$, and
set $u_n=(A\sigma)^{1/\alpha}n$. Then
$$
\im\phi^\prime(\sigma+ix\sigma)=A^{-1/\alpha}\sigma^{-\beta}x\,
\sum_{n\geq 1}\psi(u_n,x)\Delta u_n,\ref{\riemannsumref}
$$
where
$$
\psi(u,x)={1\over{2s}}{{u^\alpha\sin(xu^\alpha)}\over{\cosh u^\alpha
-\cos xu^\alpha}},\ref{\psidefref}
$$
and where $\Delta u_n=u_n-u_{n-1}=(A\sigma)^{1/\alpha}$. The sum in
{\riemannsumref} is a Riemann sum of the integral
$$
{1\over s}\int_0^\infty{{u^\alpha\sin(xe^{-u^\alpha}u^\alpha)}\over
{1+e^{-2u^\alpha}-2e^{-u^\alpha}\cos(xu^\alpha)}}\, du=C_{A,\alpha}h(x),
\ref{\integrref}
$$
where $C_{A,\alpha}=\alpha^{-1}\Gamma(\beta)\zeta(\beta)$, and where
$$
h(x)=(1+x^2)^{\beta/2}\,{{\sin(\beta\arctan x)}\over x}.\ref{\hdefref}
$$
Note that $h(x)\neq 0$, if $\alpha\geq 1$. We claim that the difference
of the Riemann sum in {\riemannsumref} and the integral {\integrref}
is $o(1)$, as $\sigma\rightarrow 0$. Indeed, this difference can be
written as
$$
\sum_{n\geq 1}\int_{u_{n-1}}^{u_n}\big(\psi(u_n,x)-\psi(u,x)\big)du,
$$
which can readily be bounded by
$$
A^{1/\alpha}\sigma^{\beta}\sum_{n\geq 1}\,\max_{u_{n-1}\leq u\leq u_n}
\,\big|{{\partial}\over{\partial u}}\psi(u,x)\big|. \ref{\intermsumref}
$$
Using the fact that, uniformly in $x$,
$$
\big|{{\partial}\over{\partial u}}\psi(u,x)\big|\leq\cases{
        O(1)u^{-1},&if $0<u\leq 1$;\cr
        O(1)e^{-(1-\epsilon)u}, &if $u>1$,}\ref{\derivestref}
$$
(with $0<\epsilon<1$), we can bound {\intermsumref} by
$$
\eqalign{
&O(1)\sigma^\beta\sum_{1\leq n\leq(A\sigma)^{-1/\alpha}}
(A\sigma)^{-1/\alpha}n^{-1}
+O(1)\sigma^\beta\sum_{n>(A\sigma)^{-1/\alpha}}e^{-(1-\epsilon)
(A\sigma)^{1/\alpha}n}\cr
&=O(1)\sigma\log (A\sigma)^{-1/\alpha}+O(1)\sigma\cr
&=o(1),\cr}
$$
and our claim follows.

\noindent
{\it Step 2}. In the general case, we write $\sigma\lambda_n=u_n
(1+\mu_n)$, with $u_n$ as before. We now claim that the difference
$$
\sum_{n\geq 1}\big(\psi(u_n(1+\mu_n),x)-\psi(u_n,x)\big)\Delta u_n
\ref{\differenceref}
$$
is $o(1)$, as $\sigma\rightarrow 0$. Indeed, using {\derivestref}
we can bound {\differenceref} by
$$
\eqalign{
&\sum_{n\geq 1}\,u_n|\mu_n|\Delta u_n\,\max_{u\in[u_{n-1},u_n]}
\,\big|{{\partial}\over{\partial u}}\psi(u,x)\big|\cr
&\leq O(1)\sigma^{1/\alpha}\sum_{1\leq\Lambda}1+O(1)\sigma^{2/\alpha}
\sum_{n>\Lambda}\mu_n ne^{-(1-\epsilon)(A\sigma)^{1/\alpha}n}\cr
&O(1)\sigma^{1/\alpha}\Lambda+O(1)\sup_{n>\Lambda}\mu_n,\cr}
$$
where $\Lambda>0$ is arbitrary. Choosing e.g. $\Lambda=
\sigma^{-1/(2\alpha)}$ we conclude that the above expression is
$o(1)$, as $\sigma\rightarrow 0$. $\square$

\vfill\eject

\centerline{\bf References}
\baselineskip=12pt
\frenchspacing

\bigskip

\item{[BC]} Bost, J.P., Connes, A.: Hecke Algebras, Type III Factors
and Phase Transitions with Spontaneous Symmetry Breaking in
Number Theory, preprint (1995)

\item{[BR]} Bratteli, O., Robinson, D.W.: {\it Operator Algebras and
Quantum Statistical Mechanics II}, Springer Verlag (1981)

\item{[C1]} Colin de Verdiere, Y.: Ergodicite et fonctions propres
du Laplacien, {\it Comm. Math. Phys.}, {\bf 102}, 497--502 (1985)

\item{[C2]} Connes, A.: {\it Noncommutative Geometry}, Academic Press
(1994)

\item{[GJ]} Glimm, J., Jaffe, A.: Boson quantum field models, in
{\it Mathematics of Contemporary Physics}, R.Streater, ed.,
Academic Press (1972)

\item{[HR]} Hewitt, E., Ross, K.A.: {\it Abstract Harmonic Analysis},
Springer Verlag (1963)

\item{[I]} Ingham, A.: A Tauberian theorem for partitions,
{\it Ann. Math.}, {\bf 42}, 1075--1090 (1941)

\item{[JLW]} Jaffe, A., Lesniewski, A., Wisniowski, M.:
Deformations of super-KMS functionals, {\it Comm. Math. Phys.}, 
{\bf 121}, 527--540 (1989)

\item{[J]} Julia, B.: Statistical Theory of Numbers, in {\it
Number Theory and Physics}, Les Houchs Winter School, J.-M. Luck,
P. Moussa and M. Waldschmidt editors, Springer Verlag (1990)

\item{[K]} Kastler, D.: Cyclic cocycles from graded KMS functionals,
{\it Comm. Math. Phys.}, {\bf 121}, 345--350 (1989)

\item{[P]} Postnikow, A.: {\it Introduction to analytic number theory},
Am. Math. Soc. (1988)

\item{[S1]} Schnirelman, A.: Ergodic properties of the eigenfunctions,
{\it Usp. Math. Nauk}, {\bf 29}, 181--182 (1974)

\item{[S2]} Spector, D.: Supersymmetry and the M\"obius inversion 
function, {\it Comm. Math. Phys.}, {\bf 127}, 239--252 (1990)

\item{[Z]} Zelditch, S.: Quantum ergodicity of $\bC^*$ dynamical
systems, {\it Comm. Math. Phys.}, {\bf 177}. 507--528 (1996)

\vfill\eject\end